\documentclass[letter,reqno,11pt] {amsart}
\usepackage{amssymb}
\usepackage{graphicx}
\usepackage[margin=1.15in]{geometry}

\newtheorem{theorem}{Theorem}[section]

\newtheorem{remark}[theorem]{Remark}
\newtheorem{prop}[theorem]{Proposition}

\newcommand{\calh}{\mathcal H}

\def\limind{\mathop{\oalign{lim\cr
\hidewidth$\longrightarrow$\hidewidth\cr}}}
\def\g{\gamma}
\def\a{\alpha}
\def\G{\Gamma}

\def\Ac{\mathcal A}
\def\Hc{\mathcal H}

\newcommand{\mysection}[1]{\section{#1}
\setcounter{equation}{0}}

\begin{document}
\title{A Survey on Rankin-Cohen Deformations}
\author{Richard Rochberg}
\address{Department of Mathematics, Washington University, St. Louis, MO, U.S.A. 63130} \email{rr@math.wustl.edu}

\author{Xiang Tang}
\address{Department of Mathematics, Washington University, St. Louis, MO, U.S.A., 63130} \email{xtang@math.wustl.edu}

\author{Yi-jun Yao}
\address{Department of Mathematics, Penn State University, State College, PA, U.S.A.,
16802} \email{yao@math.psu.edu}

\date{\today}
\maketitle
\begin{abstract}This is a survey about recent progress in
Rankin-Cohen deformations. We explain a connection between
Rankin-Cohen brackets and higher order Hankel forms.
\end{abstract}

\vspace{0.5cm}
\mysection{Introduction}

The famous \it Erlanger Programm \rm of Klein says that geometry is
about to study the transformation groups of various spaces, or more
precisely the properties invariant under the actions of such groups,
\it i.e.\rm, the symmetries.

Noncommutative geometry(NCG), which originated from Connes'
study in operator algebras in 1970's, brought the landscape of
geometry many new objects and some astonishing phenomena.

Back to the early 1990s, Connes and Moscovici pointed out that in
noncommutative geometry(NCG) while noncommutative spaces are
represented by the algebras (usually noncommutative $C^*$-algebras)
of ``continuous functions'' over noncommutative spaces, the \it
local\rm\ symmetries are reflected in some Hopf algebras. One of the
first noncommutative spaces studied in NCG is the $C^*$-algebra of a
foliated space. In the case of codimension $n$ foliations, Connes
and Moscovici discovered a Hopf algebra ${\mathcal H}_n$, which
governs the local symmetry of leaf spaces of foliations of
codimension $n$. The Hopf algebra ${\mathcal H}_n$ is universal in
the sense that it depends only on the codimension of a foliation.
This family of Hopf algebras $\{{\mathcal H}_n\}$ is very useful in
the study of transverse index theory, and later was found to have
connections with various different areas of mathematics, c.f.
\cite{CK98}, \cite{CM03-1}.  In this paper, we review the
application of the Hopf algebra ${\mathcal {H}}_1$ in Rankin-Cohen
deformations, which was initiated by Connes and Moscovici
\cite{CM03-2}.

We start by recalling the general setting of transverse geometry.
Let $M$ be a smooth manifold and $\mathcal F$ be a foliation on $M$
of codimension $n$. Let $X$ be a complete flat transversal of
$\mathcal F$, and $F^+X$ be the oriented frame bundle of $X$. The
holonomy pseudogroup $\Gamma$ acts on $X$ and therefore $F^+X$ by
transforming $X$ parallelly along paths in leaves of ${\mathcal F}$.
The ``transverse geometry" is to study the transversal $X$ along
with the action by the holonomy pseudogroup $\Gamma$.

In what follows we focus on the case when $n=1$, and define Connes-Moscovici's Hopf algebra ${\mathcal H}_1$. Now the complete
transversal $X$ is a flat 1-dim manifold; and the oriented frame bundle
$F^+X$ is diffeomorphic to $X\times {\mathbb R}^+$, and
$\Gamma$ is a discrete holonomy pseudogroup acting on $X$ as local
diffeomorphisms. We introduce coordinates $x$ on $X$ and $y$ on ${\mathbb R}^+$. The lifted action of $\Gamma$ on $F^+X$ is
\begin{equation}\label{eq:pseudogroup-action}
(x,y)\mapsto (\phi(x),\phi^\prime (x)y),\quad \phi\in\Gamma.
\end{equation}

On $F^+X$, there is a $\Gamma$-invariant volume form
$\omega=\displaystyle \frac{{\mathrm d}x\wedge {\mathrm d}y}{y^2}$,
which is also a symplectic form. This allows to consider the Hilbert space $L^2\left(\displaystyle F^+X, \frac{{\mathrm d}x\wedge
{\mathrm d}y}{y^2}\right)$ of square-integrable functions on $F^+X$. We are interested in two types of linear operators acting on this
Hilbert space $L^2\left(\displaystyle F^+X, \frac{{\mathrm d}x\wedge
{\mathrm d}y}{y^2}\right)$:

\begin{enumerate}
\item for $f\in C^\infty_c(F^+ X)$, we define $M_f: \xi\mapsto
f\xi$ for all $\xi\in L^2\left(\displaystyle F^+X, \frac{{\mathrm
d}x\wedge {\mathrm d}y}{y^2}\right)$;

\item for $\phi\in\Gamma$, we define $U_\phi:\xi\mapsto
\phi^*(\xi)=\xi\circ\phi$ for all $\xi\in L^2\left(\displaystyle
F^+X, \frac{{\mathrm d}x\wedge {\mathrm d}y}{y^2}\right)$.
\end{enumerate}

The smooth foliation algebra ${\mathcal A}_\Gamma$ is the algebra generated
by $M_f$'s and $U_\phi$'s with the relation
\[
U_\phi M_f=M_{\phi^*(f)} U_\phi.
\]
From this relation, we can say that ${\mathcal A}_\Gamma$ is the
cross product algebra $C^\infty_c(F^+ X)\rtimes \Gamma$. In
noncommutative geometry, this algebra ${\mathcal A}_\Gamma$ is
viewed as the algebra of smooth functions on the space of leaves
associated to the foliation ${\mathcal F}$ on $M$. The Hopf algebra
${\mathcal H}_1$ acts on the smooth foliation algebra by linear
operators.

By choosing a flat connection on $F^+X$, we consider two vector
fields, i) $X=y\partial_x$ as a lifting of the vector filed
$\partial_x$ on $X$, and ii) $Y=y\partial_y$ the Euler vector field
along the fiber direction. $X$ and $Y$ act on the smooth foliation
algebra ${\mathcal A}_\Gamma$ by
\[
X(f U_\phi)= X(f)U_\phi,\,\, Y(f U_\phi)= Y(f) U_\phi.
\]
It is easy to check that $Y$ is invariant under the action of
$\Gamma$, but $X$ is not:
\[
U_\phi X U_\phi^{-1} = X-y
\frac{\phi^{-1\prime\prime}(x)}{\phi^{-1\prime}(x)}Y.
\]
The failure of $X$ being $\Gamma$ invariant inspires higher operations. Define a
linear operator $\delta_1$ on ${\mathcal A}_\Gamma$ in the following
way:
\[
\delta_1(fU_\phi)= \mu_{\phi^{-1}}f U_\phi,
\]
where $\displaystyle \mu_{\phi^{-1}}(x,y)=y
\frac{\phi^{-1\prime\prime}(x)}{\phi^{-1\prime}(x)}$.

We compute $[Y,\delta_1]$, which turns out to be $\delta_1$
itself. But $[X,\delta_1]$ leads to a new operator, which we name
$\delta_2$:
\[
\delta_2(f U_\phi)= X(\mu_{\phi^{-1}})fU_\phi.
\]
Iterating the procedure of computing the commutator with $X$, we are
led to a sequence of operators $\delta_n$ acting on ${\mathcal
A}_\Gamma$ by
\[
\delta_n(f U_\phi)= X^{n-1}(\mu_{\phi^{-1}})fU_\phi.
\]

With the above preparation, we are ready to present the Hopf algebra
${\mathcal H}_1$. As an algebra, it is the universal enveloping
algebra of an infinite dimensional Lie algebra $H_1$ whose
generators are labeled as $\{X,Y,\delta_n,\ n\in \mathbb{N}\}$ with
the following commutation relations:
\[
\begin{array}{ll}
[Y, X]=X,& [X, \delta_n]=\delta_{n+1},\\

[Y,\delta_n ]
 =n\delta_n,& [\delta_n, \delta_m]=0.
\end{array}
\]

We define the following structures on $\calh_1$:
\begin{enumerate}
\item product $\cdot : \calh_1\otimes \calh_1\to \calh_1$ is defined as the
product on the universal enveloping algebra of $H_1$.
\item coproduct $\Delta:\calh_1\to \calh_1\otimes \calh_1$ is an algebra homomorphism generated by
$$
\begin{array}{l}
\Delta Y=Y\otimes 1+1\otimes Y,\\
\Delta \delta_1=\delta_1\otimes 1+1\otimes \delta_1,\\
\Delta X=X\otimes 1+1\otimes X+\delta_1\otimes Y,\\
\Delta\delta_n=[\Delta X, \Delta \delta_{n-1}].
\end{array}
$$
\item counit $\epsilon:\calh_1\to {\mathbb C}$ is defined by taking the constant component in ${\mathcal {H}}_1$.
\item antipode $S: \calh_1\to \calh_1$ is an algebra anti-homomorphism generated by
$$
S(X)=-X+\delta_1Y,\ \ S(Y)=-Y,\ \ S(\delta_1)=-\delta_1.
$$
\end{enumerate}
It is straightforward to check that $(\calh_1, \cdot, \Delta, S,
\epsilon, id)$ defines a Hopf algebra, which acts naturally on the
smooth foliation algebra ${\mathcal A}_\Gamma$.

The Hopf algebra ${\mathcal {H}}_1$ and its higher dimensional
generalizations were discovered by Connes and Moscovici \cite{CM98}
as local symmetries on ${\mathcal A}_{\Gamma}$. In the following, we
will survey some recent developments about the applications of
${\mathcal {H}}_{1}$ in the study of modular forms and in particular
Rankin-Cohen deformations.

The results reviewed in this paper are interactions between
transverse geometry and modular form theory. These two classical
``distant" subjects are mysteriously connected due to the fact that
the Hopf algebra ${\mathcal {H}}_1$ appears in both theories as the
local symmetry of the corresponding systems. What we will develop in
the last section of this paper is to add one more subject to this
story, namely the Hankel forms and transvectant theory in harmonic
analysis. There, we will introduce an algebra ${\mathcal B}_\Gamma$
associated to a pseudogroup $\Gamma$ acting on a 1-dim complex
domain by holomorphic transformations. We will show that the Hopf
algebra ${\mathcal {H}}_1$ acts on ${\mathcal B}_\Gamma$. And
interestingly, through the Hopf algebra action the Rankin-Cohen
brackets on modular form are translated to Hankel forms of higher
weights \cite{JP}.\\

\noindent{\bf Acknowledgements:} We would like to express our
gratitude to Connes and Moscovici for insightful suggestions and
constant supports in our work. We would also like to thank Pevzner,
S\"undall, and Zagier for helpful discussion.

\mysection{Rankin-Cohen brackets and deformations}

In number theory, modular forms are very important because the
coefficients of their Fourier expansions encode a great amount of
number theoretical information. We recall the definition of a
modular form. Let $\Gamma$ be a congruence subgroup of
$SL_2({\mathbb Z})$ (a subgroup of $SL_2({\mathbb R})$ such that the
entries of a matrix are all integers), a modular form of weight $2k$
is a function $f$ which satisfies:

\begin{itemize}
\item  (holomorphy) $f$ is holomorphic on the upper-half plane $\mathbb
H$,
\item (modularity) for $\gamma=\displaystyle \left(\begin{array}{cc}a & b\\
                                                                    c &
d\end{array}\right)\in
                                                                    \Gamma$
                                                                    and $z\in\mathbb
                                                                    H$,
                                                                    $f\Big|_{2k}\gamma=f$,
                                                                    where
\[
\left(f\Big|_{2k}\gamma\right)(z)=(cz+d)^{-2k}f\left(\frac{az+b}{cz+d}\right),
\]
(this means the invariance of the form $f \, dz^{k}$);
\item (growth condition at the boundary) $|f(z)|$ is assumed to be controlled by a polynomial in $\max\{1,
Im(z)^{-1}\}$.
\end{itemize}

We denote by ${\mathcal M}(\Gamma)$ the (weight) graded algebra of
modular forms with respect to the group $\Gamma$.

The Rankin-Cohen brackets are a family of universal formulas
describing all bilinear operations (up to a scalar) between spaces
of modular forms (of even weight) that can be defined in terms of
derivatives. In the 50's Rankin started the research on
bi-differential operators on ${\mathcal M}(\Gamma)$ which produce
new modular forms, and twenty years later Heri Cohen\footnote{We
have given the first name of the author because in the literature of
Rankin-Cohen deformations, there are two different and important
authors with the same last name, ``Cohen". } gave a complete answer
(cf. \cite{C75}) to this question by showing that all these
operators are linear combinations of the bracket
\begin{equation}\label{eq:rc-brackets}
[f, g]_n=\sum_{r=0}^n (-1)^r {n+2k-1 \choose n-r}{n+2l-1 \choose r}
f^{(r)}g^{(n-r)} \in {\mathcal M}_{2k+2l+2n}(\Gamma),
\end{equation}
\noindent where $f\in {\mathcal M}_{2k}$ and $g\in {\mathcal
M}_{2l}$ are two modular forms,  and
$f^{(r)}=\displaystyle\left(\frac{1}{2\pi i}\frac{\partial}{\partial
z}\right)^r f$. The above bilinear operation $[\ ,\ ]_n$ is called
$n$-th Rankin-Cohen bracket.

The original and obvious importance of the Rankin-Cohen brackets in
number theory is that these brackets often give rise to non-trivial
identities between the Fourier coefficients of modular forms, e.g.
H. Cohen's foundational paper \cite{C75}. Zagier \cite{Z94} has
found several other ``raisons d'\^etre'' for the Rankin-Cohen
brackets. He discovered the following procedure to obtain the
Rankin-Cohen brackets: we denote Ramanujan's derivation on modular
forms by $X:{\mathcal M}_{2k}(\G) \rightarrow {\mathcal
M}_{2k+2}(\G)$:
\begin{eqnarray*}
X f&=& \frac{1}{2  \pi i} \,  \frac{df}{dz} - \frac{1}{2  \pi i} \,
\frac{\partial}{\partial z} (\log \eta^4) \cdot k f,\cr \eta (z)  &
=& q^{1/24} \, \tiny{ \prod_{n=1}^{\infty} }(1 - q^n) \,  , \quad q
= e^{2\pi i z},
\end{eqnarray*}
and define two sequences of modular forms by induction: for $\Phi\in
{\mathcal M}_4(\G)$, $f_0=f, g_0=g$,
\begin{equation}
\label{eq:zagier-iteration}
\begin{split}
f_{r+1} &=X f_r + r(r+2k-1) \Phi f_{r-1},  \cr
g_{s+1} &=X g_s +
s(s+2l-1) \Phi g_{s-1}.
\end{split}
\end{equation}

Zagier showed \cite{Z94} that the Rankin-Cohen brackets can be
written as
\[
\sum_{r=0}^n (-1)^r {n+2k-1 \choose n-r}{n+2l-1 \choose r}
f_{r}g_{n-r} = [f, g]_n.
\]
These forms of the brackets have the advantage that one can easily
see the modularity of $[f,g]_n$ without any extra effort from the
definition of $f_r$ and $g_s$, while in the original presentation,
the modularity of the brackets is far from being obvious as the
derivative of a modular form is in general not a modular form any
more.

Zagier and his collaborators \cite{CMZ} showed that the collection
of all Rankin-Cohen brackets together gives rise to
(non-commutative) associative deformations of the algebra ${\mathcal
M}(\Gamma)$ of modular forms. We review their constructions in more
detail.

Paula Cohen\footnote{Here is the other author with the same last
name ``Cohen".}, Manin, and Zagier \cite{CMZ} established a
bijection between formal series with modular form coefficients and
formal invariant pseudodifferential operators. For a modular form
$f$ of weight $2k$ we define
\begin{equation}\label{D}
{\mathcal D}_{-k}(f,
\partial)=\sum_{n=0}^{\infty} \frac{\displaystyle{-k \choose
n}{-k+ -1\choose n}}{\displaystyle{-2k\choose n}}
f^{(n)}\partial^{-k-n},
\end{equation}
and the composition of two such formal pseudodifferential operators
(plus the bilinearity) give us the following product:
\begin{eqnarray}\label{tn}
f\star g&=&\sum_{n=0}^\infty c_n(k,l)[f,g]_n,\\
c_n(k,l)&=&t_n(k, l):=\frac{1}{{-2l \choose
n}}\sum_{r+s=n} \frac{{-k \choose r}{-k-1 \choose r}}{{-2k \choose
r}}\frac{{n+k+l \choose s}{n+k+l-1 \choose s}}{{2n+2k+2l-2\choose
s}}.
\end{eqnarray}
Meanwhile, by checking the first several terms, Eholzer \cite{Z94}
conjectured that
\begin{equation}\label{ehol}
f\star g :=\sum_{n=0}^\infty [f, g]_n
\end{equation}
is also an associative product.

In order to include this phenomenon in their framework, P. Cohen,
Manin, and Zagier \cite{CMZ}modified their definition of invariant
formal pseudodifferential operators, and they finally obtained a
whole family of such deformations, parametrized by the Riemann
sphere(\cite{CMZ},\cite{Z94}).   In this most general case,  we
consider $\g\in SL_2({\mathbb C})$ action on functions $\xi\in
C^\infty({\mathbb H})$:
\[
(W_\g^{{\kappa}} \xi)(z)=\xi(\frac{az+b}{cz+d})(cz+d)^{{\kappa}}.
\]
Then we can define

\begin{equation}\label{Dkappa}
{\mathcal D}_{-k}^{{\kappa}} (f,
\partial)=\sum_{n=0}^{\infty} \frac{\displaystyle{-k \choose
n}{-k+{{\kappa}} -1\choose n}}{\displaystyle{-2k\choose n}}
f^{(n)}\partial^{-k-n},
\end{equation}
which has the following invariant property:
\[
{\mathcal D}_{-k}^{{\kappa}} (f|_{2k}\g,
\partial)=W_\g^{{{\kappa}}
*}{\mathcal D}_{-k}^{{\kappa}} (f,
 \partial)W_\g^{{\kappa}}.
\]

By composing two such formal invariant pseudodifferential operators, we obtain, over the algebra of formal series with coefficients in
modular forms ${\mathcal M}[[\hbar]]$, that the natural linear
extension plus the following formula gives an associative product:
for $f\in{\mathcal M}_{2k}$, $g\in{\mathcal M}_{2l}$,
\begin{equation}\label{mukappa}
\mu^\kappa (f, g)=\sum_{n=0}^\infty t_n^\kappa(k, l) [f, g]_n
\hbar^n,
\end{equation}
where the coefficients are
\begin{equation}\label{tnkappa}
t_n^\kappa(k, l)=\frac{1}{\displaystyle{-2l\choose n}}\sum_{r+s=n}
\frac{\displaystyle{-k\choose r}{-k+\kappa-1\choose
r}}{\displaystyle{-2k\choose
r}}\frac{\displaystyle{n+k+l-\kappa\choose s}{n+k+l-1\choose
s}}{\displaystyle{2n+2k+2l-2\choose s}}.
\end{equation}
Moreover, the coefficients $t_n^\kappa(k, l)$ are conjectured to be
equal to
\begin{equation}\label{tnkappa1}
t_n^\kappa(k,l)=\left(-\frac{1}{4}\right)^n \sum_{j\geq 0} {n\choose
2j}\frac{ \displaystyle{-\frac{1}{2}\choose
j}{\kappa-\frac{3}{2}\choose j}{\frac{1}{2}-\kappa\choose j}}
{\displaystyle{-k-\frac{1}{2}\choose j}{-l-\frac{1}{2}\choose
j}{n+k+l-\frac{3}{2}\choose j}}.
\end{equation}

By taking $\kappa=\displaystyle\frac{1}{2}$ or
$\displaystyle\frac{3}{2}$ under this form, (1.3) turns to be
Eholzer's product Eq. (\ref{tn}), and until 2004 this is the only
possible way to prove the associativity of Eholzer's product.
(Unfortunately, Zagier's original proof \cite{Z} of the identity Eq.
(\ref{tnkappa1}) is not published, and one can find an elementary
but rather long proof in \cite{these}).

\mysection{Modular Hecke algebras and Connes-Moscovici's deformation}

The interaction between the theory of modular forms and
noncommutative geometry goes back to December 2001, when Zagier gave
a course at Coll\`ege de France(\cite{Z01}) and Connes was in the
audience. One year later, Connes and Moscovici discovered that the
Hopf algebra ${\mathcal H}_1$ that controls the local symmetry of
the transverse geometry of codimension one foliations does act on
some big algebra constructed from the algebra of modular forms, and
named it ``modular Hecke algebra", which we will briefly recall now.

We first define
\begin{equation*}
{\mathcal M} (\Gamma(N)) := \Sigma^{\oplus} {\mathcal M}_{2k}
(\Gamma(N)) \,  ,  \, {\mathcal M}^{0} (\Gamma(N)) :=
\Sigma^{\oplus} {\mathcal M}^{0}_{2k} (\Gamma(N)),
\end{equation*}
These algebras form a projective system with respect to the divisibility of the integer $N$. Define
\begin{equation*}
{\mathcal M} :=  \limind_{N \rightarrow \infty} \,  {\mathcal M}
(\Gamma(N)) \, ,  \quad \hbox{resp.} \quad {\mathcal M}^{0} :=
\limind_{N \rightarrow \infty} \,  {\mathcal M}^0 (\G(N)) \,  .
\end{equation*}

An {\it operator Hecke form of level $\Gamma$} \cite{CM03-1} is a map
\begin{eqnarray*}
 F: \G\backslash GL_2^+ ({\mathbb Q}) & \rightarrow &  {\mathcal M} \,  , \\
 \Gamma_\alpha & \mapsto & F_{\alpha} \in {\mathcal M}  \,  ,
\end{eqnarray*}
which has a finite support, and satisfies the {\it covariance
condition}:
\begin{equation*}
F_{\alpha \g} \,  (z) = \,  F_{\alpha} \vert \g(z) , \quad \forall \alpha \in GL_2^+ ({\mathbb Q})
, \gamma \in \G ,  z\in {\mathbb H} \,  .
\end{equation*}
The modular Hecke algebra  ${\mathcal A}(\Gamma)$ is an associative
algebra consisting of operator Hecke forms of level $\Gamma$ with
the product,
\begin{equation*}
(F^1 * F^2)_\alpha : =\sum_{\beta\in \Gamma\backslash
GL_2^+({\mathbb Q})} F^1_\beta\cdot
F^2_{\alpha\beta^{-1}}\Big|\beta.
\end{equation*}

An important discovery of Connes and Moscovici in \cite{CM03-1} is
that the Hopf algebra ${\mathcal H}_1$ acts on ${\mathcal
A}(\Gamma)$.

Before we give the detail of this action, we briefly recall the
general definition of a Hopf algebra action on an algebra. Let $H$
be a Hopf algebra and $M$ be an algebra. We say that
$\alpha:H\otimes M\to M$ defines an action of $H$ on $M$ if the
following two conditions hold;
\begin{enumerate}
\item $M$ is an $H$-module with respect to the algebra structure on $H$;
\item the following property holds with respect to the coalgebra structure on $H$,
\[
\alpha(h, a_1a_2)=m\big((\alpha\otimes \alpha) (\Delta(h), a_1\otimes a_2)\big),
\]
where $a_1,a_2$ are elements of $M$, $\Delta$ is the coproduct of $H$, and $m: M\otimes M\to M$ is the multiplication operator.
\end{enumerate}
With this definition, one can easily check that Connes-Moscovici's
Hopf algebra $\calh_1$ acts on the smooth foliation algebra
${\mathcal {A}}_\Gamma$.  In the following, we review the
$\calh_1$'s action on modular Hecke algebra ${\mathcal {A}}(\Gamma)$
as was introduced by Connes-Moscovici \cite{CM03-1}.

For an $f\in {\mathcal M}_{2k}$, we define
\begin{equation*}
Xf = \frac{1}{2  \pi i} \frac{\partial}{\partial z} - \frac{1}{2 \pi
i} \, \frac{\partial}{\partial z} (\log \eta^4) \cdot kf\, ,\,Y(f) =
k \cdot f .
\end{equation*}
For  $ \g = \left(\begin{array}{cc} a &b\\
 c &d  \end{array}\right) \in GL^+ (2, \mathbb{Q})\ $, we define
\begin{equation*} \label{nmu}
\begin{split}
 \mu_{\g} \,  (z) &= \frac{1}{2 \,  \pi^2}
 \left(  G^*_2 | \g \,  ( z) -  G^*_2 ( z)
+ \frac{2  \pi i \,  c}{cz+d} \right),\\
\hspace{-1cm}G_2^* (z) &= 2\zeta (z)+ 2 \sum_{m\geq 1}\sum_{n\in
{\mathbb Z}}\frac{1}{(mz+n)^2} =\frac{\pi^2}{3}- 8\pi^2 \sum_{m,
n\geq 1}m e^{2\pi i m n z}.
\end{split}
\end{equation*}
We notice that here $\,  \mu_{\g}\equiv 0 \, $ if $\,  \g \in
SL_2 (\mathbb{Z}) \, $.

We define the action of $X,Y, \delta_n$ on ${\mathcal {A}}(\Gamma)$ as follows, for $ F \in \Ac (\G) \,  ,  \a \in G^+ (\mathbb{Q})$,
\begin{eqnarray*}\label{actn}
X(F)_{\a} \,  &:=& \,  X(F_{\a}) \,  ,   \qquad   \cr
 Y(F)_{\a} \, &:=& \,  Y(F_{\a}) \,  ,  \qquad  \cr
 \delta_n (F)_{\a}  &:=& \,  \mu_{n,  \a} \cdot F_{\a} \,  ,\ \text{where\ }\mu_{n,  \,  \a} \,  := \,  X^{n-1}(\mu_{\a}) \,  ,
 \quad  n \in {\mathbb N}.
\end{eqnarray*}
With the above preparation, it is not difficult to check the
following theorem.
\begin{theorem}\label{thm:hopf-action-modular}(\cite{CM03-1})
Let $\G$ be a congruence subgroup of $SL_2(\mathbb{Z})$.
\begin{enumerate}
\item The Hopf algebra
 $\Hc_1$ acts on the algebra $\Ac (\G)$.
\item The Schwarzian derivative $\displaystyle \delta'_2=\delta_2
-\frac{1}{2}\delta_1^2$ is inner and is implemented by
$\displaystyle \omega_4=-\frac{1}{72} E_4 \in \Ac (\G) \,  $, where
\[
E_4=1+240\sum_{n=1}^\infty n^3 \frac{e^{2\pi i nz}}{1-e^{2\pi i nz
}}
\]
is an Eisenstein series of degree 4.
\end{enumerate}
\end{theorem}
In fact, Connes and Moscovici \cite{CM03-1} pointed out that
${\mathcal {A}}(\Gamma)$ can actually be obtained as the crossed
product of the algebra of modular forms by the Hecke ring, and then
reduced by a projection which is determined by the congruence
subgroup $\Gamma$. Based on this fact, and assuming the
associativity of the Eholzer product,  Connes and Moscovici
\cite{CM03-2} subsequently proved that this associative formal
deformations can be canonically extended from the algebra of modular
forms to the modular Hecke algebra associated to a congruence
subgroup.

\begin{remark}
Inspired by Theorem \ref{thm:hopf-action-modular} (2),  Connes and
Moscovici  introduced a concept of projectivity of an ${\mathcal
H}_1$ action on an algebra $\mathcal A$ as follows: there is an
element $\Omega \in {\mathcal A}$ such that,
\begin{enumerate}
\item
\begin{equation} \label{om}
\delta'_2(a):=\left(\delta_2-\frac{1}{2}\delta_1^2\right)(a)= \Omega
\, a - a \, \Omega\, , \:\forall a \in \Ac;
\end{equation}
\item due to the commutativity of the
$\delta_k$'s,
\begin{equation} \label{om1}
   \delta_k(\Omega )=\, 0\, , \: \forall k \in \mathbb{N}.
\end{equation}
\end{enumerate}
We point out that this ${\mathcal {H}}_1$ projectivity structure is
a generalization of the projective structure on elliptic curves. On
an elliptic curve, a projective structure means a choice of atlas
such that the transition function between different charts can be
chosen to be in $SL_2(\mathbb{R})$. The Ramanujan differential
appears as a connection associated to such a projective structure.
\end{remark}

A crucial observation of Connes and Moscovici \cite{CM03-2} is that
the extended Rankin-Cohen brackets on the modular Hecke algebra
${\mathcal A}(\Gamma)$ can be represented using elements in
${\mathcal {H}}_1$ and the element $\omega_4\in {\mathcal
A}(\Gamma)$.  For example, the first Rankin-Cohen bracket can be
realized by
\[
RC_1(a\otimes b)=m((X\otimes Y-Y\otimes X-\delta_1Y\otimes Y)(a\otimes b)),
\]
where $m: {\mathcal A}(\Gamma)\otimes {\mathcal A}(\Gamma)\to
{\mathcal A}(\Gamma)$ is the multiplication on ${\mathcal
A}(\Gamma)$. Generalizing the above formula of $RC_1$, we define
\begin{eqnarray*}
\hspace{-0.5cm}A_{n+1}&:=& S(X)\,  A_n - n\, \Omega^{\, o} \,
\left(Y- \frac{n-1}{2}\right) A_{n-1},\cr \hspace{-0.5cm}B_{n+1}&:=&
X\,  B_n - n\, \Omega \, \left(Y- \frac{n-1}{2}\right) B_{n-1},
\end{eqnarray*}
where $A_{-1}:=0,  \,  A_{0}:=1$ and $B_0:=1$,  $B_1:= X$, and
$\Omega^o$ is the multiplication operator from the right by
$\Omega$. We notice that this is very similar to the sequences Eq.
(\ref{eq:zagier-iteration}) of modular forms that Zagier used to
reformulate the classical Rankin-Cohen brackets with
$\Omega=\omega_4$. In general, the $n$-th Rankin-Cohen bracket can
be written as
 \begin{equation*}\label{eq:cm-def}
RC_n(a, b):= \sum_{k=0}^n \,  \frac{A_k }{k!} \, (2Y+k)_{n-k}(a)\,
\,  \frac{ B_{n-k}}{(n-k)!} \,  (2Y+n-k)_{k}(b).
\end{equation*}

We remark that the above $RC_n$ is a generalization of the classical
$n$-th Rankin-Cohen bracket in the sense that when we take the
modular Hecke algebra and we restrict the $n$-th Rankin-Cohen
bracket on modular forms, we get the classical $n$-th Rankin-Cohen
bracket. The advantage of this generalization is that now we can
apply $RC_n$ to an arbitrary algebra $A$ on which ${\mathcal {H}}_1$
acts with a projective structure.

By a technique called \it full injectivity\rm\ , Connes and
Moscovici \cite{CM03-2} proved their main theorem which states that
the associative product can ultimately be lifted to a universal
deformation formula for projective actions of the Hopf algebra
${\mathcal H}_1$ :

\begin{theorem}(\cite{CM03-2})\label{thm:univ-def-cm}  The functor $RC_*:=\sum RC_n$
applied to any algebra $\Ac$ endowed with a projective structure
yields a family of formal associative deformations of $\Ac$, whose
products are given by
\[
f\star g =\sum RC_n(f, g) \hbar^n.
\]
\end{theorem}

The full injectivity method used in the proof of the above theorem
essentially says that there are enough different actions of
${\mathcal H}_1$ with projective structures so that any cocycle properties which lead to an associative deformation (on the algebras on which ${\mathcal H}_1$ acts)
can always be lifted to the Hopf algebra level.

\mysection{Rankin-Cohen deformation via Fedosov}
                        \label{Fedosov}

In the previous section, we have seen a beautiful result of Connes
and Moscovici extending structures in modular form theory to study
the Hopf algebra ${\mathcal {H}}_1$. In this section, we look at the
universal deformation formula obtained in Theorem
\ref{thm:univ-def-cm} from the view point of transverse geometry of
codimension one foliation.

In geometry, deformation of the algebra of smooth functions on a
manifold has been studied for a long time. In particular, it is not
difficult to see that the first order limit of an associative
deformation of the commutative algebra of smooth functions on a
manifold defines a Poisson bracket on the manifold. Here, by a
Poisson bracket on a manifold $P$ we mean a bilinear map $\{\ ,\ \}:
C^\infty(P)\otimes C^\infty(P)\to C^\infty(P)$ such that for any
$f,g,h\in C^\infty(P)$,
\begin{enumerate}
\item $\{f,g\}=-\{g,f\}$,
\item $\{f, gh\}=g\{f,h\}+\{f,g\}h$,
\item $\{\{f,g\},h\}+\{\{h,f\},g\}+\{\{g,h\},f\}=0.$
\end{enumerate}
When $P$ has a symplectic structure, namely, a non-degenerate closed
 2-form $\omega$ on $P$, then the inverse of $\omega$ defines a
Poisson bracket on $C^\infty(P)$ by $\{f,g\}:= \omega^{-1}(df,dg)$.
($\omega$ is viewed as a skew symmetric bilinear form on $TM$, and
its inverse defines a skew symmetric bilinear form on $T^*M$.) In
mathematical physics, the phase space of a physical system is
usually described by a symplectic manifold, and observables of a
physical system are smooth functions on the symplectic manifold.
Bayer-Flato-Fronsdal-Lichnerowicz-Sternheimer \cite{BFFLS} pointed
out that we can use deformations of the algebra of smooth functions
on a symplectic (Poisson) manifold, the phase space of a physical
system, to study the corresponding quantum system. They call a
deformation of the algebra $C^\infty(P)$ with the first order term
equal to the Poisson structure a {\it deformation quantization} of
the symplectic (Poisson) manifold. An easy and beautiful example of
such a theory is that in quantizing a free particle on $\mathbb{R}$,
we have the algebra of quantum observables are generated by the
position operator $\hat{q}$ (the multiplication operator by function
$q$) and the momentum operator $i\hbar \frac{d}{dx}$ on
$L^2({\mathbb {R}})$. Such an algebra is a deformation of the
algebra of smooth functions on $\mathbb{R}^2$ (the corresponding
phase space) with the standard symplectic structure $dp\wedge dq$.
The product of this deformation quantization of
$C^\infty(\mathbb{R}^2)$ can be written as
\begin{equation}\label{eq:moyal}
f\star
g(x)=\exp\left(-\frac{i\hbar}{2}\omega^{ij}\frac{\partial}{\partial
y^i}\frac{\partial}{\partial z^j}\right)f(y)g(z)|_{x=y=z},\qquad
f,g\in C^\infty(\mathbb {R}^2),
\end{equation}
which is called the Moyal product \cite{Moy49}. With this geometry
and physics in mind, we ask ourselves whether there is a geometric
(maybe physical) interpretation of the Rankin-Cohen deformation.

We start with looking for a better understanding of a simplified
version of the Rankin-Cohen deformation obtained in Theorem
\ref{thm:univ-def-cm}. We notice that if we set all $\delta_i$
($i=1,\cdots, \infty$) and $\Omega$ (the projective structure) to be
zero in the Rankin-Cohen bracket $RC_n$, we obtain a Universal
Deformation Formula(UDF) of the 2-dimensional solvable Lie algebra
$h_1$ associated to the $ax+b$ group, i.e. $h_1=\text{span}\{X,Y\}$
with $[Y,X]=X$. We call the simplified deformation the {\it reduced
Rankin-Cohen product}:
\begin{equation}\label{RCreduit}
\hspace{-1.7cm}RC_{red} = \sum_{n=0}^{\infty}\frac{\hbar^n}{n!}
\sum_{k=0}^n\left[ (-1)^k{n\choose k}X^k(2Y+k)_{n-k}\otimes
X^{n-k}(2Y+n-k)_{k} \right],
\end{equation}
with $Y_k=Y(Y+1)\cdots (Y+k-1)$.

The $ax+b$ acts on the upper half place by translation,
$(x,y)\mapsto (ax+b, ay)$. This induces an action of $h_1$ on
$C^\infty(\mathbb{R}\times \mathbb{R}^+)$. Now applying the reduced
Rankin-Cohen bracket (\ref{RCreduit}) on $C^\infty(\mathbb{R}\times
\mathbb{R}^+)$, we obtain a deformation of the algebra
$C^\infty(\mathbb{R}\times \mathbb{R}^+)$. By using an argument via
the orbit method (more precisely by a theorem of Gutt \cite{G83}),
we together with Bieliavsky \cite{BTY} proved that the reduced
Rankin-Cohen deformation on $C^\infty(\mathbb{R}\times
\mathbb{R}^+)$ is isomorphic to the Moyal product (\ref{eq:moyal})
over the half-plane.

In mid 90's, Giaquinto and Zhang proposed another UDF for
$h_1$ (cf. \cite{GZ98}):
\begin{equation*}\label{GZ}
F=\sum_{n=0}^{\infty}\frac{\hbar^n}{n!}\sum _{r=0}^n (-1)^r
{n\choose r}X^{n-r}Y_{r}\otimes X^{r}Y_{n-r}.
\end{equation*}
With the same idea, we \cite{BTY} were able to prove that the
Gianquinto-Zhang UDF applied to $C^\infty(\mathbb{R}\times
\mathbb{R}^+)$ is also isomorphic to Moyal product as well. We
remark that the action of the universal enveloping algebra of $h_1$
on $C^\infty(\mathbb{R}\times \mathbb{R}^+)$ is fully injective.
Therefore, the associativity of the Moyal product can be used to
prove that the realization of the reduced Rankin-Cohen deformation
on $C^\infty(\mathbb{R}\times \mathbb{R}^+)$ is associative. By the
full injectivity assumption, we can give a new proof that the
reduced Rankin-Cohen deformation and also the Giaquinto-Zhang
deformation are associative.

The above study provides a connection between the reduced
Rankin-Cohen brackets and the well-known Moyal product. This also
inspires us to pursue further whether we can give a geometric
reconstruction of Connes-Moscovici's Rankin-Cohen deformation of the Hopf algebra ${\mathcal {H}}_1$.

When asking ourselves this question, we realized that there is
already a potential direction for the answer. Remember ${\mathcal
{H}}_1$ was introduced in the study of transverse geometry of
codimension one foliation. In particular, ${\mathcal {H}}_1$ acts on
the smooth foliation algebra ${\mathcal {A}}_\Gamma$. In
noncommutative geometry, ${\mathcal {A}}_\Gamma$ is viewed as the
algebra of smooth functions on the quotient space $F^+X/\Gamma$,
which is usually a non-hausdorff topological space. A key
observation here is that the volume form $dx\wedge dy/y^2$ defines a
symplectic form on $F^+X$, which is invariant under the action of
$\Gamma$. This leads us to the consideration of deformation
quantization of $F^+X$ with the symmetry of the pseudogroup
$\Gamma$.

As a first step, we need to understand the first order term of the
Rankin-Cohen deformation. As we have mentioned, when we consider the
deformation of $C^\infty(P)$, the first order term is a Poisson
structure. Now, for a noncommutative algebra ${\mathcal
{A}}_\Gamma$, what is the first order term of a deformation? The
answer to this question is that, it is a noncommutative Poisson
structure.

A noncommutative Poisson structure on an algebra $A$ is a degree 2
Hochschild cocycle $\Pi$ such that the Gerstenhaber bracket $[\Pi,
\Pi]_G$ is a coboundary.

We are able to prove the following result by a long but direct
computation.
\begin{theorem}\label{thm:poisson}(\cite{BTY}) Let $A$ be an algebra equipped with an ${\mathcal {H}}_1$ action.
Then $\Pi(a,b)=X(a)Y(b)-Y(a)X(b)+\delta_1(Y(a))Y(b)=m(RC_1(a\otimes
b))$ defines a noncommutative Poisson structure on $A$.
\end{theorem}

In the second step, we aim to understand the geometric meaning of a
projective structure. As we mentioned, classically a projective atlas on a riemann surface
$X$ assign a principal $SL_2(\mathbb{R})$ bundle on $X$. In the case
of codimension 1 foliation, we obtained the following theorem giving
a geometric interpretation of a projective structure on ${\mathcal
A}_\Gamma$ analogous to this classical picture.
\begin{theorem}\label{thm:connection}(\cite{BTY})
In the standard action of ${\mathcal {H}}_1$ on a smooth foliation
algebra ${\mathcal {A}}_\Gamma$, the projective structure is
equivalent to the existence of an invariant connection on $F^+X$ of
the form
\[
\begin{split}
\nabla_{\partial_x}\partial_x=\mu(x,y)\partial_y,\qquad\qquad&
\nabla_{\partial_x}\partial_y=\frac{1}{2y}\partial_x,\\
\nabla_{\partial_y}\partial_x=\frac{1}{2y}\partial_x,\qquad\qquad&
\nabla_{\partial_y}\partial_y=-\frac{1}{2y}\partial_y.
\end{split}
\]
\end{theorem}

We point out that the connection introduced in Theorem
\ref{thm:connection} is a symplectic torsion free connection on
$F^+X$. This reminds us Fedosov's theory \cite{fed} of deformation
quantization of a symplectic manifold.

We briefly explain Fedosov's construction of a deformation
quantization on the manifold $(F^+X, \omega)$. We consider a Weyl
algebra bundle $W$ whose fiber at every point $p$ of $F^+X$ is the
Weyl algebra $W_p$ consisting of formal power series
\[
a(u, \hbar)=\sum\limits_{k, |\alpha|\geq 0} \hbar^k a_{k, \alpha}
u^\alpha.
\]
Here $\hbar$ is the formal parameter, $y=(u^1, u^{2})\in T_x F^+X$,
and $\alpha=(\alpha_1, \alpha_{2})\in {\mathbb Z}_{\geq 0}\times
{\mathbb Z}_{\geq 0}$, and
$y^\alpha=(u^1)^{\alpha_1}(u^{2})^{\alpha_{2}}$.

The product on $W_p$  is the Moyal product, for $a,b\in W_p$
\begin{eqnarray*}
\hspace{-1cm} a\circ b&=&
\exp\left(-\frac{i\hbar}{2}\omega^{ij}\frac{\partial}{\partial
v^i}\frac{\partial}{\partial w^j}\right) a(v, \hbar)b(w,
\hbar)|_{v=w=u}\cr  &=& \sum_{k=0}^{\infty}
\left(-\frac{i\hbar}{2}\right)^k \frac{1}{k!}\omega^{i_1
j_1}\cdots\omega^{i_k j_k}\frac{\partial^k a}{\partial u^{i_1}\cdots
\partial u^{i_k}}\frac{\partial^k b}{\partial u^{j_1}\cdots \partial
u^{j_k}}.
\end{eqnarray*}

Let $\bigwedge$ be the bundle $\wedge^\bullet T^* F^+X$. An abelian
connection $D: \Gamma^{\infty}(W\otimes \bigwedge)\to
\Gamma^{\infty}(W\otimes \bigwedge)$ is a connection on $W$ with $
D^2 a=0$, for any $a\in \Gamma^\infty(W\otimes \bigwedge)$. When
considering $W_D:=\{a, Da=0\}$, Fedosov proved that

\noindent{\bf Theorem.}  (\cite{fed}){\it \label{thm:fedosov}
For all $a_0\in C^\infty(F^+X)[[\hbar]]$, there exists a unique
section $a\in W_D$, noted as $\sigma^{-1}(a_0)$, such that
$\sigma(a):=a(x, 0, \hbar)=a_0$.  Hence, $\sigma$ is a bijection
between $W_D$ and $C^{\infty}(F^+X)[[\hbar]]$. And we can define on
$ C^{\infty}(F^+X)[[\hbar]]$ an associative product}
\[
\label{eq:star} a \star b=\sigma(\sigma^{-1}(a)\circ
\sigma^{-1}(b)).
\]

In the case that $F^+X$ has a $\Gamma$-invariant symplectic
connection, the second named author proved \cite{T04} that the above
deformation quantization of $C^\infty_c(F^+X)$ defines a deformation
of the corresponding crossed product algebra ${\mathcal
A}_\Gamma=C^\infty_c(F^+X)\rtimes \Gamma$.

When we are given the $\Gamma$-invariant symplectic torsion free
connection $\nabla$ (Thm. \ref{thm:connection}) on $F^+X$, if we
assume furthermore that $\mu(x,y)$ is of the form $y\nu(x)$ with
$\nu(x)$ an arbitrary smooth function of variable $x$, the
connection is actually flat. This allows to find an explicit formula
for the abelian connection $D$. The equation $Da=0$ gives us a
system of differential equations. And starting from $a_{0, 0}=f$ and
solving the system by induction, we obtain two sequences of
elements:
\begin{eqnarray}
 A_{m+1}&=-X A_m -
m\displaystyle\frac{\mu}{y^3}\left(Y-\frac{m-1}{2}\right)A_{m-1}, \\
\displaystyle  B_{m+1}&=X B_m -
m\displaystyle\frac{\mu}{y^3}\left(Y-\frac{m-1}{2}\right)B_{m-1}.
\end{eqnarray}
where $\displaystyle  X=\frac{1}{y}\frac{\partial}{\partial
x}\,\,,\,\, Y=-y\frac{\partial}{\partial y}$. Extending the above
formulas onto the smooth foliation algebra ${\mathcal A}_\Gamma$, we
obtained the same recurrence relation that Connes and Moscovici used
in their definition of generalized Rankin-Cohen brackets, so one can
consider what we get as another realization of the Rankin-Cohen
deformation. The advantage of this argument is that one gets the
associativity without any extra effort. We show \cite{BTY} that in
this context we can also prove the ``full injectivity", which allows
us to conclude the associativity of the Rankin-Cohen deformation at
the Hopf algebra level.

We notice that in Thm \ref{thm:poisson}, the first Rankin-Cohen
bracket always defines a Poisson structure no matter whether there
is an ${\mathcal {H}}_1$ projective structure or not. This inspires
to ask whether ${\mathcal {H}}_1$ has a Universal Deformation
Formula without assuming the existence of a projective structure. In
the language of deformation quantization, this question is whether
there is a deformation of the algebra ${\mathcal A}_\Gamma$ without
the existence of an invariant symplectic connection. We learned from
Fedosov \cite{fed} and also Gorokhovsky-Bressler-Nest-Tsygan
\cite{BFFLS} that the existence of an invariant symplectic torsion
free connection is not a necessary condition for a smooth foliation
algebra to have a deformation quantization. Actually, we can always
construct a deformation quantization of a smooth foliation algebra
${\mathcal A}_\Gamma$ using the idea of algebroid stacks. With this
in mind, we proved the following theorem, where we can drop the
assumption of the existence of a projective structure.
\begin{theorem}(\cite{TY})
The Hopf algebra $\calh_1$ has a universal deformation formula, i.e.
there is an element $R\in {\mathcal H}_1[[\hbar]]\otimes_{{\mathbb
C}[[\hbar]]} {\mathcal H}_1[[\hbar]]$ satisfying
\begin{eqnarray}
((\Delta\otimes 1)R) (R\otimes 1)&=&((1\otimes \Delta)R) (1\otimes
R),\\
(\epsilon\otimes 1)(R)&=&1\otimes 1=(1\otimes \epsilon) (R).
\end{eqnarray}
\end{theorem}
As a side remark, we also obtained \cite{TY} a proof of the associativity of the Eholzer product of a reasonable length using some elementary methods.

\mysection{Rankin-Cohen via Representations}\label{rep}

Rankin-Cohen brackets and related deformation questions can also be
studied using the theory of infinite dimensional representations of
$SL_2(\mathbb{R})$. First of all, one can explicitly give an
interpretation of these brackets using unitary representation theory
of $SL_2({\mathbb R})$. This fact is known by experts for a long
time\footnote{For example, Deligne wrote in 1973 (cf. \cite{D73}) :
Remarque 2.1.4. L'espace $F(G, GL_2({\mathbb Z}))$ ci-dessus est
stable par produit. D'autre part, $D_{k-1}\otimes D_{l-1}$ contient
les $D_{k+l+2m-1} (m \geq 0)$ . Pour $ m = 0$ , ceci correspond au
fait que le produit $fg$ d'une forme modulaire holomorphe de poids
$k$ par une de poids $l$, en est une de poids $k+l$ . Pour $m = 1$ ,
en coordonn\'ees (1.1.5.2), on trouve que $\displaystyle l
\frac{\partial f}{\partial z}.g - kf.\frac{\partial g}{\partial z}$
est modulaire holomorphe de poids $k+l+2$ , et ainsi de suite. De
m\^eme dans le cadre ad\'elique. } but does not seem to be clearly
written anywhere. The main result is the following theorem:

\noindent {\bf Theorem}. {\it Let $f\in{\mathcal M}_{2k},
g\in{\mathcal M}_{2l}$ be two modular forms. Let  $\pi_f\cong
\pi_{\deg f}, \pi_g\cong \pi_{\deg g}$ be associate representations
which are discrete series of the group $SL_2({\mathbb R})$. The
tensor product of these two representations is decomposed into a
direct sum of discrete series,

\[
\pi_f\otimes \pi_g =\bigoplus_{n=0} \pi_{\deg f+\deg g+2n}.
\]

\noindent The Rankin-Cohen bracket $[f, g]_n$ gives (up to a scalar)
the minimal $K$-weight vectors in the representation space of the
component $\pi_{\deg f+\deg g+2n}$.}

We remark that the result about tensor product between $SL_2(
{\mathbb{R}})$ representations are well-known (c.f. \cite{Rep78} by
Repka). The new part in the above theorem is the relation between
Rankin-Cohen brackets and representation theory.

The representation $\pi_f$ is constructed in the following way: let
$f\in {\mathcal M}_{2k}(\Gamma)$ be a modular form, one can
associated to it a function over $\Gamma\backslash SL_2({\mathbb
R})$ by using the following mapping:

\[
(\sigma_{2k} f)(\g)=f_{|_k}|\g(i)=(ci+d)^{-2k}
f\left(\frac{ai+b}{ci+d}\right), \, \, \, \,  \text{for}\, \,
\g=\left(\begin{array}{cr}
                           a & b\\
                            c & d
                            \end{array}\right)\in SL(2,{\mathbb R}).
\]
This function belongs to
$$
C^\infty(\Gamma\backslash SL_2({\mathbb R}), 2k)=\{F\in
C^\infty(\Gamma\backslash SL_2({\mathbb R})), F(\g\cdot r_\theta)=
\exp( i 2k\theta) F(\g) \},
$$
where
\[
r_\theta=\left(\begin{array}{cc}\cos(\theta)&-\sin(\theta)\\
\sin(\theta)&\cos(\theta)\end{array}\right).
\]

By taking into account the natural action of $SL_2( {\mathbb R})$ on
$C^\infty(\Gamma\backslash SL_2({\mathbb R}))$:
\[
(\pi(h) F)(g)=F(gh),
\]
one gets a representation of $SL_2({\mathbb R})$ and also a
representation of the complexified Lie algebra ${\mathfrak
sl}_2({\mathbb C})$ by taking the smallest invariant subspace which
contains the orbit of $\sigma_{2k} f$. One can show that this
representation is a discrete series of weight $2k$. In this step, we
take all the vectors in a base of the representation space to a
subspace of $C^\infty( {\mathbb H})$, by using the inverse of
$\sigma_{2(k+n)}, n\geq 0$.

\smallskip

With this interpretation, we studied the deformed products in a more
general setting. From now on we drop the holomorphy condition on a
modular form, and consider all functions which satisfy the
modularity condition with respect to $\Gamma$. We denote the algebra
under consideration by $\widetilde{M}(\Gamma)$. Enlarging
$\widetilde{M}(\Gamma)$, we consider the tensor algebra
\[
\widetilde{\mathcal M}(\Gamma)^\otimes =\sum_n \widetilde{\mathcal
M}(\Gamma)^{\otimes n},
\]
and define
\begin{eqnarray}\label{scmdef}
{\textsc m} : \widetilde{\mathcal M}(\Gamma)^\otimes &\rightarrow &
\widetilde{\mathcal M}(\Gamma), \cr f_1\otimes f_2\otimes\cdots \otimes
f_n &\mapsto & f_1f_2\cdots f_n,
\end{eqnarray}
and we extend the degree operator by
\[
\deg(f_1\otimes f_2\otimes\cdots \otimes f_n )=\sum_{i=1}^n \deg
f_i.
\]

We define a derivation on $\widetilde{{\mathcal M}}(\Gamma)$, for an element $f\in \widetilde{\mathcal M}(\Gamma)$ of degree $2k$,
\begin{equation}\label{deftildepar}
\tilde{X} f =\frac{1}{2\pi i}\frac{df}{dz} - \frac{2kf}{4\pi Im(z)}.
\end{equation}
We extend $\tilde{X}$ to $\widetilde{\mathcal M}(\Gamma)^\otimes$ by the Leibnitz rule
\[
\tilde{X} (f_1\otimes f_2\otimes\cdots \otimes f_n )=\sum_{i=1}^n
f_1\otimes f_2\otimes\cdots \otimes\tilde{X}f_i\otimes\cdots \otimes
f_n.
\]

Inspired by the Rankin-Cohen brackets (\ref{eq:rc-brackets}), we study the following two families of formal products. The first family is $\star :
\widetilde{\mathcal M}(\Gamma)^\otimes [[\hbar]]\otimes \widetilde{\mathcal
M}(\Gamma)^\otimes [[\hbar]]\rightarrow \widetilde{\mathcal
M}(\Gamma)^\otimes [[\hbar]]$ defined by linear extension and the
following formula:
\begin{equation}\label{prodform}
f\star g =\sum \frac{A_n(\deg f,  \deg g)}{(\deg f)_n(\deg g)_n}
\left(\sum_{r=0}^n (-1)^r \tilde{X}^r {\deg f+n-1\choose n-r} f
\otimes \tilde{X}^{n-r} {\deg g+n-1\choose r}g\right) \hbar^n,
\end{equation}

\noindent where $f, g\in \widetilde{\mathcal M}(\Gamma)^\otimes$.
The second family of products $\ast$ is the restriction of $\star$
defined by Eq. (\ref{prodform}) to $\widetilde{\mathcal
M}(\Gamma)\subset \tilde{\mathcal M}(\Gamma)^\otimes$ composed with
the application $\textsc m$ defined in (\ref{scmdef}). More
concretely,  $\ast: \widetilde{\mathcal M}(\Gamma) [[\hbar]]\times
\widetilde{\mathcal M}(\Gamma) [[\hbar]]\rightarrow
\widetilde{\mathcal M}(\Gamma) [[\hbar]]$ is defined as follows:
\begin{eqnarray}\label{prodformfaib}
f\ast g & = & {\textsc m}(f\star g)\cr &=&\sum \frac{A_n(\deg f,
\deg g)}{(\deg f)_n(\deg g)_n} \left(\sum_{r=0}^n (-1)^r \tilde{X}^r
{2k+n-1\choose n-r} f \tilde{X}^{n-r} {2l+n-1\choose r}g\right)
\hbar^n\cr &=& \sum \frac{A_n(\deg f,  \deg g)}{(\deg f)_n (\deg
g)_n} [f, g]_n \hbar^n,
\end{eqnarray}
where $f, g\in \widetilde{\mathcal M}$, and the notation
$(\alpha)_n:=\alpha(\alpha+1)\cdots(\alpha+n-1)$. We will always take the (natural) assumption $A_0=1$ and $A_1(x,y)=xy$, and $A_i$ $(i\geq 2)$ is a polynomial of two variables.

Our main goal is to study the associativity of the above two family
of products. We name the associativity of the product $\star$ {\it
strong associativity}, and the associativity of $\ast$ {\it weak
associativity}.

The third named author proved \cite{Y} that, while it is obvious
that a family of coefficient functions $A_n(\deg f,\deg g)$ which
defines a strongly associative product also defines also a weakly
associative product, the inverse, which is not clear {\it a priori},
is also true. This equivalence between the two types of
associativity is used to show the following results,
\begin{theorem}(\cite{Y}) Cohen-Manin-Zagier \cite{CMZ} have found {\it all} associative formal deformation $\ast:\widetilde{\mathcal
M}[[\hbar]]\times\widetilde{\mathcal M}[[\hbar]]\rightarrow
\widetilde{\mathcal M}[[\hbar]]$ defined by the linearity and the
formula
\begin{eqnarray}
f\ast g & = & \sum \frac{A_n(\deg f,  \deg g)}{(\deg f)_n (\deg
g)_n} [f, g]_n \hbar^n,
\end{eqnarray}
\noindent where $\widetilde{\mathcal M}(\Gamma)$ is the space of the
functions which satisfy the modularity condition, and the notation
$(\alpha)_n:=\alpha(\alpha+1)\cdots(\alpha+n-1)$. One assumes again
$A_0=1$ and $A_1(x, y)=xy$.
\end{theorem}

The main reason for this claim to be valid is that the determinant
of some $2\times 2$ linear system is non-zero, and the drop of
holomorphy allows us to have the total freedom to modify a function
in the interior of a fundamental domain. Furthermore, with the help
of some computation of certain multivariable polynomials done by
Mathematica, the third named author was able to prove the following
proposition.
\begin{prop}(\cite{Y}) Let $\Gamma$ be a congruence
subgroup of  $SL_2({\mathbb Z})$ such that ${\mathcal M}(\Gamma)$
admits the unique factorization property (for example $SL_2({\mathbb
Z})$ itself), let $F_1, F_2, G_1, G_2\in {\mathcal M}(\Gamma)$ such
that
\begin{equation*}
RC(F_1, G_1)=RC(F_2, G_2),
\end{equation*}
as formal series in ${\mathcal M}(\Gamma)[[\hbar]]$, then
there exists a constant $C$ such that
\begin{equation*}
F_1=C F_2,  G_2= C G_1.
\end{equation*}
\end{prop}

\mysection{Hankel forms and a new ${\mathcal {H}}_1$ action}
\label{hankel}

In this section, we aim to set up a connection between the Hopf
algebra ${\mathcal {H}}_1$ and the theory of high order Hankel forms
introduced by Janson and Peetre\footnote{In literature, there is a
huge amount of study of Hankel forms and transvectants. Due to our
limited knowledge, we have only cited those references related to
our work.} \cite{JP}. We hope that such a connection will inspire
more interactions among transverse geometry, number theory, and
harmonic analysis. Similar to what we did in Sec. \ref{rep}, we will
consider a set of functions more general than modular forms. Instead
of dropping the holomorphy property of a modular form, in this
section we will drop the modularity property. For example, we will
consider holomorphic functions on the unit disk of ${\mathbb {C}}$
that are square integrable with respect to some measure.

Let $D$ be the unit disk in the complex plane $\mathbb {C}$.  Let
$\eta(z)=1-|z|^2$, and $d\sigma(z)=(1/\pi)dxdy$. The weighted
Bergman space $A^{2\alpha}(D)$ ($\alpha>-1$) is defined to be
\[
\{f: \bar{\partial}f=0, \int_{D}|f(z)|^2\eta^\alpha d\sigma(z)<\infty\}.
\]
The group $SL(2, \mathbb{R})$ acts on the weighted Bergman space
$A^{2\alpha}(D)$ as follows, for
$\gamma=\left(\begin{array}{cc}a&b\\ c&d\end{array}\right)$,
\begin{equation}\label{eq:action-weighted}
\pi_{2k}(\gamma)(f)(z)=f(\frac{az+b}{cz+d})(cz+d)^{-\alpha}.
\end{equation}

It is easy to check that above action $\pi_{2k}$ defines a unitary
representation of $SL(2, \mathbb{R})$, which is actually
irreducible. Such a representation is called a discrete series of
$SL(2, \mathbb{R})$. Applying the result of Repka \cite{Rep78}, a
tensor product of representations in the discrete series is a direct
sum of irreducible $SL(2, \mathbb{R})$ representations in the
discrete series. More explicitly, as $SL(2, {\mathbb R})$
representations, we have the following decomposition,
\[
\pi_{2k}\otimes \pi_{2l}=\bigoplus_{i\geq 0}\pi_{2k+2l+2i}.
\]

We denote the projection from $A^{2k}(D)\otimes A^{2l}(D)$ to
$A^{2k+2l+2i}$ by $\Pi_{i}$, which is $SL_2({\mathbb R})$
equivariant bilinear map from $A^{2k}(D)\times A^{2l}(D)$ to
$A^{2k+2l+2i}$. Such an equivariant bilinear form is called by
Janson and Peetre \cite{JP} a Hankel form of weight $2i$. Similar to
the Rankin-Cohen brackets, $P_{2i}$ can be expressed using the
derivatives and the weights of the two components, i.e. up to a
constant,
\[
\Pi_{2i}(f,g)(z)=\sum_{s=0}^{i}(-1)^{i-s}\left(\begin{array}{l}i\\
s\end{array}\right)\frac{1}{(w(f))_{s}(w(g))_{i-s}}\partial^sf\partial^{i-s}g,
\]
where $w(f)$ and $w(g)$ are weights of $f$ and $g$.

We remark that the above discussion of Hankel forms of high weights
is closely related to the theory of transvectant. Olver-Sanders
\cite{OS}, Choi-Mourrain-Sol\'e \cite{CMS}, and El Gradechi
\cite{e06} provided an explicit relation between transvectants,
Rankin-Cohen brackets, and Moyal product \cite{OS}. One can even
uses these interesting relations to give a proof of the
associativity of the Rankin-Cohen deformation, which was explained
by Pevzner \cite{P}. We refer readers to the cited references for
more details.

In the following of this section, we aim to introduce an algebra
${\mathcal B}_\Gamma$ associated to a pseudogroup $\Gamma$ of
holomorphic transformations on a complex domain $\Sigma$ of
dimension one. And we will show that there is a natural action of
${\mathcal {H}}_1$ on ${\mathcal B}_\Gamma$. In the special case
when $\Gamma$ is trivial and $\Sigma$ is $D$, the unit disk in
${\mathbb C}$, the Rankin-Cohen deformation gives rise to Hankel
forms of higher weights.

Let $\Sigma$ be a complex domain of dimension 1 (maybe of several
components). $\Gamma$ is a pseudogroup acting on $\Sigma$ by
holomorphic transformations. Let $T^{1,0}\Sigma$ be the holomorphic
tangent bundle of $\Sigma$, which is equipped with a $\Gamma$ action
\[
\gamma(x,y)=(\gamma(x), \partial_x\gamma(x)y),
\]
where $(x,y)$ is a point on $\sigma$. We assume that the holomorphic
tangent bundle is trivialized. Namely, we have chosen a global
coordinates $(x,y)$ on $T^{1,0}\Sigma$.

\begin{remark}
We have required the bundle $T^{1,0}\Sigma$ to be trivialized in
order to work with holomorphic functions and holomorphic vector
fields on $\Sigma$. If $T^{1,0}\Sigma$ is not a trivial bundle over
$\Sigma$, then we will have to work with sheaves of holomorphic
functions and holomorphic vector fields. Accordingly, we will have
to enlarge the notion of Hopf algebra to Hopf algebroid like in
\cite{CM01}.
\end{remark}

We consider the space of holomorphic functions on $T^{1,0}\Sigma$,
and denote it by $A_{T^{1,0}\Sigma}$. The action of $\Gamma$ on
$T^{1,0}\Sigma$ induces an action of $\Gamma$ on
$A_{T^{1,0}\Sigma}$. Therefore, we define ${\mathcal B}_\Gamma$ to
be the crossed product algebra $A_{T^{1,0}\Sigma}\rtimes \Gamma$.

We remark that when $\Sigma$ is the unit disk in ${\mathbb C}$, the
space of holomorphic functions on $T^{1,0}D$ has a natural
decomposition,
\[
A_{T^{1,0}D}=\bigoplus_{i\geq 0}A_Dy^i,
\]
where $y$ is the coordinate along the fiber direction of $T^{1,0}D$,
and $A_D$ is the space of holomorphic functions on $D$.
$SL_2(\mathbb{R})$, the group of holomorphic transformations of $D$,
acts on each component $A_Dy^i$ exactly like the $SL_2( \mathbb{R})$
action on the weighted Bergman space $A^{2i}(D)$. Actually both
$A_D$ and $A^{2i}(D)$ share a dense subalgebra, the algebra of
polynomials of variable $x$.

To introduce the Hopf algebra ${\mathcal {H}}_1$ action on
$B_\Gamma$, we first introduce two holomorphic vector fields on
$T^{1,0}\Sigma$. Define
\[
X=y\partial_x,\qquad Y=y\partial_y.
\]
It is easy to check that for an element $\gamma\in \Gamma$,
\[
U_\gamma X
U_{\gamma^{-1}}=X-y\frac{\partial_x^2\gamma^{-1}}{\partial_x\gamma^{-1}}Y,\qquad
U_\gamma Y U_{\gamma^{-1}}=Y,
\]
where $U_\gamma$ is the action of $\gamma$ on $A_{T^{1,0}\Sigma}$.

The vector fields $X$ and $Y$ lift to act on the algebra ${\mathcal
B}_\Gamma$ by acting on the component of $A_{T^{1,0}\Sigma}$. We
define $\delta_1:{\mathcal B}_\Gamma\to {\mathcal B}_\Gamma$ by
\[
\delta_1(fU_\gamma)=y\partial_x\big(\log(\partial_x\gamma^{-1})\big)fU_\gamma.
\]
It is not difficult to check that operators $X,Y,\delta_1$ satisfy
\[
[Y,X]=X,\quad[Y,\delta_1]=\delta_1,
\]
and $[X, \delta_1]=\delta_2$, where $\delta_2:{\mathcal B}_\Gamma\to
{\mathcal B}_\Gamma$ is an operator defined by
\[
\delta_2(fU_\gamma)=y^2\partial_x^2\big(\log(\partial_x(\gamma^{-1}))\big).
\]
Continuing this procedure, we have $\delta_n:{\mathcal B}_\Gamma\to
{\mathcal B}_\Gamma$ by
\[
\delta_n(fU_\gamma)=X\big(\delta_{n-1}(fU_\gamma)\big).
\]
\begin{prop}\label{prop:hol}Connes-Moscovici's Hopf algebra ${\mathcal
{H}}_1$ acts naturally on the algebra ${\mathcal B}_\Gamma$.
\end{prop}
\begin{proof}
The proof that ${\mathcal {H}}_1$ acts on ${\mathcal B}_\Gamma$ is a
repetition of Connes-Moscovici's proof that ${\mathcal {H}}_1$ acts
on the smooth foliation algebra ${\mathcal {A}}_\Gamma$. Here we are
replacing smooth functions by holomorphic functions and
differentiations by holomorphic differentiations. But the Hopf
algebraic structure involved does not change at all.
\end{proof}

As an application of the Rankin-Cohen deformation, assuming $\Gamma$
is trivial, we apply the reduced Rankin-Cohen deformation
(\ref{RCreduit}) to the algebra $A_{T^{1,0}D}$, we have for
$f(x)y^k,g(x)y^l$,
\[
\begin{split}
RC_{red}&(f(x)y^k,
g(x)y^l)=\sum_{n=0}^\infty\frac{\hbar^n}{n!}\\
&\sum_{s=0}^n(-1)^s\left[\left(\begin{array}{c}n\\
s\end{array}\right)\frac{(2k+n-1)!}{(2k+s-1)!}\partial_x^{s}(f)(x)\frac{(2l+n-1)!}{(2l+n-s-1)!}\partial_x^{n-s}(g)(x)\right]y^{k+l+n}.
\end{split}
\]
When we take the component of $\hbar^n$, then we find
\[
\begin{split}
R_n(fy^k,
gy^l)&=(-1)^n\frac{(2k+n-1)!(2l+n-1)}{((2k-1)!(2l-1)!)}\\
&\sum_{s=0}^n(-1)^{n-s}\left[\left(\begin{array}{c}n\\
s\end{array}\right)\frac{1}{(2k)_s(2l)_{n-s}}\partial_x^{s}(f)(x)\partial_x^{n-s}(g)(x)\right]y^{k+l+n}.
\end{split}
\]
Recalling the relation between $A_{T^{1,0}D}$ and the weighted
Bergman space $A^{2n}(D)$, we can extend the action of $h_1$ onto
$\oplus_{n\geq 0}A^{2n}(D)$ and conclude with the following
proposition.
\begin{prop}\label{prop:hankel}
The reduced $i$-th Reduced-Cohen bracket on $\oplus_{n\geq
0}A^{2n}(D)$ is a Hankel form of weight $2i$. The associativity of
the reduced Rankin-Cohen deformation implies these family of Hankel
forms defines an associative deformation of $\oplus_{n\geq
0}A^{2n}(D)$.
\end{prop}

We remark that in \cite{UU} (and also \cite{dp}),
Unterberger-Unterberger introduced some operator symbol calculus
associated to weighted Bergman spaces. And they are also able to
construct an associative deformation on $\oplus_{n\geq 0}A^{2n}(D)$.
It is an interesting question to compare our deformation in
Proposition \ref{prop:hankel} and the deformation obtained in
\cite{UU}, while our coefficients is much simpler than the ones in
\cite{UU}.

\end{document}